\documentclass[12pt]{amsart}

\usepackage[centertags]{amsmath}
\usepackage{amsfonts}
\usepackage{amsthm}
\usepackage{newlfont}
\usepackage{amscd}
\usepackage{amsgen}
\usepackage{amssymb}
\usepackage{stmaryrd}
\usepackage{amssymb,amsmath}
\usepackage{amscd}
\usepackage{enumerate}

\newlength{\defbaselineskip} 
\theoremstyle{plain}
\newtheorem{thm}{Theorem}[section]
\newtheorem{cor}[thm]{Corollary}
\newtheorem{con}[thm]{Conjecture}
\newtheorem{df}[thm]{Definition}
\newtheorem{lema}[thm]{Lemma}
\newtheorem{obs}[thm]{Proposition}
\newtheorem{exm}[thm]{Example}

\newtheorem{rem}[thm]{Remark}

\theoremstyle{definition} 
\theoremstyle{definition} \newtheorem{ex}{Example}[section] %

 \numberwithin{equation}{section}
\def\p{\mathbb{P}}

\def\n{\mathbb{N}}

\def\p{\mathbb{P}}

\def\ob{\begin{obs}}
\def\kob{\end{obs}}
\def\dow{\begin{proof}}
\def\kdow{\end{proof}}

\def\tw{\begin{thm}}
\def\ktw{\end{thm}}
\def\hip{\begin{con}}
\def\khip{\end{con}}
\def\lem{\begin{lema}}
\def\klem{\end{lema}}
\def\ex{\begin{exm}}
\def\prog{\begin{pr}}
\def\kprog{\end{pr}}
\def\wn{\begin{cor}}
\def\kwn{\end{cor}}
\def\uwa{\begin{rem}}
\def\kuwa{\end{rem}}
\def\kex{\end{exm}}
\def\dfi{\begin{df}}
\def\kdfi{\end{df}}

\begin{document}
\title{Family of counterexamples to King's conjecture}
\author{Mateusz Micha\l ek}
\maketitle
\begin{abstract}
In this short note we present an infinite family of arbitrary high dimensional counterexamples to the King's conjecture.
\end{abstract}
\section{Introduction}
In \cite{king} King made the following conjecture:

\hip For any smooth, complete toric variety $X$ there
exists a full, strongly, exceptional collection of line bundles.
\khip
The conjecture turned out to be false. The first counterexample was given by
Hille and Perling, in \cite{Hille-Perling}. They showed
a smooth, complete toric surface that does not have a full
strongly exceptional collection of line bundles. Although very nice techniques were used, the proof was burdensome and required large tables of cases. Later counterexamples of higher dimension were found and the case of toric surfaces was classified \cite{nowy}. The main purpose of this note is to present a rather compact, combinatorial proof that an infinite family of higher dimensional varieties does not satisfy King's conjecture. The varieties we consider were suspected to be counterexamples and these are $\p^n$ blown up in two points. It turns out that for $n$ sufficiently large the longest strongly exceptional collection of line bundles is (by a factor) shorter then the rank of the Grothendieck group, hence such a collection is far from being full.

For more information on this topic the reader is advised to look in \cite{rosa-costa-tilting} and references therein.
\section*{Acknowledgements}
The author would like to thank very much Rosa Maria Mir\'o-Roig and Laura Costa for inspiring talks, examples and ideas. I also thank Markus Perling for important remarks and careful reading of the article. 
\section{Notation}
The varieties we consider are of Picard number 3. Such (smooth, complete, toric) varieties have been fully classified by Batyrev in \cite{baty} in terms of their {\it primitive
collections}. Using this classification our case is $|X_0|=|X_2|=|X_3|=|X_4|=1$ and $|X_1|=n-1$ with all other parameters equal to $0$.
To each element of $X_i$ corresponds a $T$-divisor. All divisors in a given $X_i$ are linearly equivalent and are given by $D_v$, $D_y$, $D_z$, $D_t$, $D_u$ respectively for $i=0,1,2,3,4$.

In the strongly exceptional collection differences of divisors cannot have nonzero higher cohomology. Divisors with nonzero higher cohomology will be called forbidden.
The following classification of forbidden divisors is very easy to establish. In a general case of Picard number three this has been done in \cite{ML}, but in this special case one can use arguments of elementary topology.
The forbidden divisors in our case are $\alpha_1 D_v+\alpha_2 D_y+\alpha_3 D_z+\alpha_4 D_t+\alpha_5 D_u$, where exactly $2,3$ or $5$ consecutive (in a cyclic way, that is indices are considered modulo 5) $\alpha$'s are negative and if $\alpha_2<0$, then $\alpha_2\leq -n+1$.

We have $D_z=D_t+D_y$ and $D_v=D_u+D_y$. We choose the basis $D_y$, $D_t$, $D_u$, what gives us forbidden divisors $(\alpha_1+\alpha_2+\alpha_3)D_y+(\alpha_3+\alpha_4)D_t+(\alpha_1+\alpha_5)D_u$ with the conditions on $\alpha$'s as above. A divisor $aD_y+bD_t+cD_u$ will be denoted by $(a,b,c)$ and we reserve precise letters for precise coordinates.
A line bundle $L_1$ will be called compatible with $L_2$ iff they can both appear in a strongly exceptional collection, that is iff $L_1-L_2$ and $L_2-L_1=-(L_1-L_2)$ are not forbidden.
\section{Proof}
Let us fix a strongly exceptional collection $E$. We assume without loss of generality that $0\in E$ and that all other divisors in $E$ have nonnegative coefficient $a$.
\lem\label{l1}
The only divisors with $a=0$ compatible with $(0,0,0)$ are: $$(0,-1,0),(0,0,-1),(0,1,0),(0,0,1),(0,-1,1),(0,1,-1).$$
\klem
\dow
If $b<-1$, then we take $\alpha_1=0$, $\alpha_2=1$, $\alpha_3=-1$ $\alpha_4$-negative to obtain $b$, $\alpha_5$-any to obtain $c$. Analogously for $c<-1$,  hence $-1\leq b,c\leq 1$. Moreover $(0,-1,-1)$ is also bad (so also $(0,1,1)$).
\kdow
\wn\label{trzy}
There can be at most 3 distinct line bundles with $a=0$ in $E$. 
For a fixed $a$ we can have only 3 line bundles in $E$.
\kwn
\dow
Follows by inspection.
\kdow
\lem\label{gl}
For $a>0$ the only line bundles $(a,b,c)$ that are not forbidden must satisfy $-1\leq b\leq a$ and $-1+a-b\leq c\leq a$ (and by symmetry $-1 + a - c \leq b \leq a$).
\klem
\dow
For $b<-1$ we take $\alpha_1=0$, $\alpha_3=-1$, $\alpha_2$-positive to have $a$, $\alpha_4$-negative to have $b$, $\alpha_5$-any to have $c$. For $b>a$ we look at $(-a,-b,-c)$ and take $\alpha_3=-a$, $\alpha_1=\alpha_2=0$, $\alpha_4$-negative to have $-b$, $\alpha_5$-any to have $-c$. In the same way $-1\leq c\leq a$\footnote{The parameters $b$ and $c$ are in symmetry.}. So the only case that we have to exclude is $-1\leq c<-1+a-b$. In such a case we can take $\alpha_4=-1$, $\alpha_3=b+1$, $\alpha_2=0$, $\alpha_1=a-b-1$, $\alpha_5=c-a+b+1<0$.
\kdow
\lem\label{8}
For three consecutive parameters $a$'s there can be at most $8$ line bundles in $E$.
\klem
\dow
We assume without loss of generality $0\leq a\leq 2$. If the lemma does not hold, then from \ref{trzy} we would have to have 3 line bundles for each $a$.
For $a=0$ we can have either:

Case 1: $(0,0,0),(0,-1,0),(0,0,-1)$ then for $a=1$ there is only one compatible from \ref{gl} namely $(1,0,0)$.

Case 2:$(0,0,0),(0,1,0),(0,0,1)$ then for $a=1$ the compatible line bundles are $(1,1,1),(1,1,0),(1,0,1)$. If we choose all of them then the only one compatible for $a=2$ is $(2,1,1)$ from \ref{gl}.

Case 3: $(0,0,0),(0,-1,0),(0,-1,1)$; $(0,0,0),(0,0,-1),(0,1,-1)$;\linebreak
$(0,0,0),(0,1,0),(0,1,-1)$; $(0,0,0),(0,0,1),(0,-1,1)$. All these possibilities are cases 1 or 2 after subtracting a divisor from all three considered divisors.
\kdow
\dfi
Line bundles in the collection $E$ with $a>n$ are called high. Others are called low.
\kdfi
\lem\label{jeden}
A high line bundle is forbidden unless either $b=1$ (high bundles of type 1)or $c=1$ (high bundles of type 2).
\klem
\dow
Suppose that $b=0$ or $b=-1$. We show that $(-a,-b,-c)$ is forbidden. Take $\alpha_1=-1$, $\alpha_2=-a+1$, $\alpha_3=0$, $\alpha_4=-b$, $\alpha_5$-any to obtain $-c$. So $b\geq 1$ and analogously $c\geq 1$. If both coefficients were strictly greater then $1$ we would obtain $(-a,-b,-c)$ by taking all $\alpha$'s negative.
\kdow
\lem\label{jedwE}
We cannot have high line bundles of both types in $E$.
\klem
\dow
From \ref{gl} a high line bundle must have the coordinate different from $1$ greater or equal to $n-1$. If we subtract two high line bundles of different types we can assume that the first coordinate is positive and one of the others will be less or equal to $-n+2$ what contradicts \ref{gl} for $n>3$.
\kdow
From now on without loss of generality we assume that we only have high line bundles of type 1 in $E$.
Let us project all high line bundles from $E$ onto the first coordinate obtaining a subset of $\n$. Suppose that this subset has got $k$ elements, that is high line bundles can have $k$ different parameters $a$. We obtain:
\lem\label{bound}
There are at most $k+2$ high line bundles in $E$.
\klem
\dow
We assumed that $0\in E$, so the high line bundles in $E$ must not be forbidden. We know that for each high line bundle in $E$ we have $b=1$, so from
\ref{gl} we know that $0\leq a-c\leq 2$. Let us notice that the difference $a-c$ cannot decrease when $a$ increases for high line bundles in $E$. Indeed suppose that we have two high line bundles in $E$ of the form $(a_1,1,c_1)$, $(a_2,1,c_2)$ with $a_2>a_1$ and $a_2-c_2<a_1-c_1$. By subtracting these two line bundles we obtain $(a_2-a_1,0,c_2-c_1)$ that is forbidden by \ref{gl}.

Notice that each time we have more then one line bundle for a fixed $a$ then the difference $a-c$ strictly increases. This means that we can have one line bundle for each $a$ plus possibly two more as $a-c$ increases from $0$ to $2$. This gives us in total $k+2$ line bundles.

\kdow
\ob
There are at most $\frac{8}{3}(n-1)+6$ low line bundles (from \ref{8}), so $k>0$ for $n>13$.
\kob
\uwa
Of course $k$ is at most $n+1$. Otherwise we would have two high line bundles in $E$ with the difference that is high. By \ref{jedwE} the difference would have $b=0$, hence by \ref{jeden} it would have $c=1$ and would be forbidden by \ref{gl}.
\kuwa
From the definition of $k$ we know that there is a line bundle $L=(a,1,c)$ in $E$, with $a\geq n+k$. Now we investigate line bundles with $a<k$, that are called very low.
\lem\label{pom}
Each very low line bundle in $E$ must have $b=0$.
\klem
\dow
Let $B$ be a very low line bundle. $L-B$ is high, so from \ref{jeden} either the second or third coordinate is $1$. The third one is $c_L-c_B\geq a_L-2-a_B> n+k-2-k=n-2>1$, for $n>3$. We see that $b_L-b_B=1$. As $b_L=1$ the theorem follows.
\kdow
For very low line bundles in $E$ the parameter $c$ is either $a$ or $a-1$ by \ref{gl} and \ref{pom}. Reasoning analogously to \ref{bound}, we see that there are at most $k+1$ very low line bundles (the difference $a-c$ cannot decrease).
\tw
The sequence $E$ can have at most: $k+1+\frac{8}{3}(n-k-1)+6+k+2\leq\frac{8}{3}n-\frac{2}{3}k+\frac{19}{3}<3n-1$ for $n>20$.
\ktw
\uwa
The bounds on $n$ can be easily improved. For example by considering separately the case $k=1$ one can decrease the bound to $n>18$. We concentrated rather on brevity of the proof than sharp bounds.
\kuwa

Mateusz Micha\l ek

{Mathematical Institute of the Polish Academy of Sciences,

\'{S}w. Tomasza 30, 31-027 Krak\'{o}w, Poland}
\vskip 2pt
{Institut Fourier, Universite Joseph Fourier,

100 rue des Maths, BP 74, 38402 St Martin d'H\`eres, France}

e-mail address:\emph{wajcha2@poczta.onet.pl}
\end{document}